\RequirePackage{ifpdf}
\ifpdf 
\documentclass[pdftex]{sigma}
\else
\documentclass{sigma}
\fi

\begin{document}

\allowdisplaybreaks

\renewcommand{\thefootnote}{$\star$}

\renewcommand{\PaperNumber}{092}

\FirstPageHeading

\ShortArticleName{Sonine Transform Associated to
the Dunkl Kernel on the Real Line}

\ArticleName{Sonine Transform Associated to
the Dunkl Kernel\\ on the Real Line\footnote{This paper is a contribution to the Special
Issue on Dunkl Operators and Related Topics. The full collection
is available at
\href{http://www.emis.de/journals/SIGMA/Dunkl_operators.html}{http://www.emis.de/journals/SIGMA/Dunkl\_{}operators.html}}}

\Author{Fethi SOLTANI}

\AuthorNameForHeading{F.~Soltani}

\Address{Department of Mathematics, Faculty of Sciences of Tunis,\\ Tunis-El Manar University, 2092 Tunis,
Tunisia}
\Email{\href{mailto:Fethi.Soltani@fst.rnu.tn}{Fethi.Soltani@fst.rnu.tn}}

\ArticleDates{Received June 19, 2008, in f\/inal form December 19,
2008; Published online December 26, 2008}

\Abstract{We consider the Dunkl intertwining operator $V_\alpha$ and its dual ${}^tV_\alpha$, we def\/ine and study the Dunkl Sonine
 operator and its dual on $\mathbb{R}$. Next, we introduce complex powers of the Dunkl Laplacian $\Delta_\alpha$ and
establish inversion formulas for the Dunkl Sonine opera\-tor~$S_{\alpha,\beta}$ and its dual ${}^tS_{\alpha,\beta}$. Also, we
give a Plancherel formula for the operator~${}^tS_{\alpha,\beta}$.}


\Classification{43A62;  43A15; 43A32}

\section{Introduction}

In this paper, we consider the Dunkl operator $\Lambda_\alpha$,
$\alpha > -1/2$,  associated  with  the ref\/lection group~$\mathbb{Z}_2$ on $\mathbb{R}$. The operators were in general
dimension introduced by Dunkl in \cite{Dunkl1} in connection with a
generalization of the classical theory of spherical harmonics; they
play a major role in various f\/ields of mathematics \cite{Dunkl2,Dunkl3,deJeu} and also
in physical applications~\cite{Lapointe}.

The Dunkl analysis with respect to $\alpha \geq -1/2$ concerns the
Dunkl operator $\Lambda_\alpha$, the Dunkl transform ${\cal
F}_\alpha$ and the Dunkl convolution $\ast_{\alpha}$ on
$\mathbb{R}$. In the limit case $(\alpha = -1/2)$;
$\Lambda_\alpha$, ${\cal F}_\alpha$ and $\ast_{\alpha}$ agree with
the operator $d/dx\,$, the Fourier transform and the standard
convolution respectively.

 First, we study the Dunkl Sonine operator $S_{\alpha,\beta}$, $\beta>\alpha$:
 \[S_{\alpha, \beta}(f)(x)
:=\frac{\Gamma(\beta+1)}{\Gamma(\beta-\alpha)\Gamma(\alpha+1)}\int^{1}_{-1}f(xt)(1
- t^2)^{\beta-\alpha - 1}(1 + t)|t|^{2\alpha+1}dt,\]  and its
 dual ${}^tS_{\alpha,\beta}$ connected with these operators. Next, we establish for them the
same results as those given in \cite{Ludwig,Solmon} for the Radon transform and
its dual; and in \cite{Nessibi} for the spherical mean operator and its dual on
$\mathbb{R}$. Especially:
\begin{enumerate}\itemsep=0pt
\item[--] We def\/ine and study the complex powers
for the Dunkl Laplacian $\Delta_\alpha=\Lambda^2_\alpha$.

\item[--] We give
inversion formulas for $S_{\alpha,\beta}$ and ${}^tS_{\alpha,\beta}$
associated with integro-dif\/ferential and
integro-dif\/ferential-dif\/ference operators when applied to some
Lizorkin spaces of functions (see \cite{Nessibi,Rachdi,Samko}).
\item[--] We establish a
Plancherel formula for the operator ${}^tS_{\alpha,\beta}$.
\end{enumerate}

The content of this work is the following.
 In Section~\ref{section2}, we recall some results about the Dunkl operators. In particular, we
give some properties of the operators $S_{\alpha,\beta}$ and
${}^tS_{\alpha,\beta}$.

 In Section~\ref{section3}, we consider the
tempered distribution $|x|^\lambda$ for $\lambda \in
\mathbb{C}\backslash \{-(\ell+1),\;\ell \in \mathbb{N}\}$ def\/ined by
\[\langle|x|^{\lambda}, \varphi \rangle : =
\int_{\mathbb{R}}|x|^{\lambda}\varphi(x)dx.\] Also we study the complex
powers of the Dunkl Laplacian $(-\Delta_\alpha)^\lambda$, for some
complex number~$\lambda$. In the classical case when $\alpha=-1/2$,
the complex powers of the usual Laplacian are given in~\cite{Stein}.

 In Section~\ref{section4}, we give the following inversion formulas:
 \begin{gather*}
 g =
S_{\alpha,\beta} K_1({}^tS_{\alpha,\beta})(g),\qquad f =
({}^tS_{\alpha,\beta})K_2S_{\alpha,\beta}(f),
\end{gather*}
where
\begin{gather*}
K_1(f)
=\frac{c_\beta}{c_\alpha}\,(-\Delta_\alpha)^{\beta-\alpha}f,\qquad
K_2(f)=\frac{c_\beta}{c_\alpha}\,(-\Delta_\beta)^{\beta-\alpha}f\qquad \mbox{and}\qquad c_\alpha=\frac{1}{[2^{\alpha+1}\Gamma(\alpha+1)]^2}.
\end{gather*}
Next, we give the following Plancherel formula for the operator
${}^tS_{\alpha,\beta}$:
\begin{gather*}
\int_{\mathbb{R}} |f(x)|^2 |x|^{2\beta+1}dx
= \int_{\mathbb{R}} |K_3({}^tS_{\alpha,\beta}(f))(y)|^2
|x|^{2\alpha+1}dy,
\end{gather*}
where
\begin{gather*}
 K_3(f)=\sqrt{\frac{c_\beta}{c_\alpha}}\,
(-\Delta_\alpha)^{(\beta-\alpha)/2}f.
\end{gather*}

\section{The Dunkl
intertwining operator and its dual}\label{section2}

 We consider the Dunkl operator $\Lambda_\alpha$,
$\alpha \geq -1/2$,  associated  with  the ref\/lection group
$\mathbb{Z}_2$ on $\mathbb{R}$:
\begin{equation}
\Lambda_\alpha f(x):= \frac{d}{dx}
f(x) + \frac{2\alpha + 1}{x} \left[ \frac{f(x) -
f(-x)}{2}\right].\label{eq1}
\end{equation}

For $\alpha \geq -1/2$ and $\lambda \in \mathbb{C}$, the initial
problem: \[
\Lambda_\alpha f(x) = \lambda f(x),\qquad  f(0) = 1,
\]
has a unique analytic solution $E_\alpha(\lambda x)$ called Dunkl kernel
\cite{Dunkl2,deJeu} given by \[
E_\alpha(\lambda x) = \Im_\alpha(\lambda x) +
\frac {\lambda x}{2(\alpha + 1)} \Im_{\alpha + 1}(\lambda x),
\]
where
 \[\Im_\alpha(\lambda x) : = \Gamma(\alpha
+ 1) \sum^{ \infty}_{n=0} \frac{(\lambda x/2)^{2n}}{n!\,\Gamma(n +
\alpha + 1)},\]
is the modif\/ied spherical Bessel function of order
$\alpha$.

 Notice that in the case $\alpha=-1/2$, we have
\[\Lambda_{-1/2} = d/dx \qquad \mbox{and} \qquad  E_{-1/2} (\lambda x) = e^{\lambda x}.\]

For $\lambda \in \mathbb{C}$ and $x\in \mathbb{R}$, the Dunkl kernel
$E_\alpha$ has the following Bochner-type represen\-ta\-tion (see
\cite{Dunkl2,Rosler1}):
\[E_\alpha(\lambda x)=a_\alpha
\int^{1}_{-1}e^{\lambda x t}\big(1 - t^2\big)^{\alpha - 1/2}(1 +
t)dt,
\]
where
\[
a_\alpha=\frac{\Gamma(\alpha+1)}{\sqrt{\pi}\,\Gamma(\alpha+1/2)},\]
which can be written as:
\begin{gather*}
E_\alpha(\lambda x) = a_\alpha\,
\mbox{sgn}(x)\,|x|^{-(2\alpha+1)}\int^{|x|}_{-|x|}e^{\lambda y}\big(x^2
- y^2\big)^{\alpha - 1/2}(x + y)dy,\qquad x\neq 0,\\
E_\alpha(0)=1.
\end{gather*}
We notice that, the Dunkl kernel $E_\alpha(\lambda x)$ can be also expanded
in a power series \cite{Rosenblum} in the form:
\begin{equation}
E_\alpha(\lambda x) = \sum^{
\infty}_{n=0} \displaystyle{\frac{(\lambda
x)^n}{b_n(\alpha)}},\label{eq2}
\end{equation}
where
\[b_{2n}(\alpha) = \frac{2^{2n}
n!}{\Gamma(\alpha+1)} \Gamma(n + \alpha + 1),\qquad b_{2n+1}(\alpha)
= 2(\alpha+1)b_{2n}(\alpha+1).\]

Let $\alpha >-1/2$ and we def\/ine the Dunkl intertwining operator
$V_\alpha$ on ${\cal E}(\mathbb{R})$ (the space of
$C^\infty$-functions on $\mathbb{R}$), by
\[V_\alpha(f)(x):=a_\alpha
\int^{1}_{-1}f(xt)\big(1 - t^2\big)^{\alpha - 1/2}(1 + t)dt,\] which can be
written as:
\begin{gather*}
V_\alpha(f)(x)= a_\alpha\,
\mbox{sgn}(x)\,|x|^{-(2\alpha+1)}\int^{|x|}_{-|x|}f(y)\big(x^2 -
y^2\big)^{\alpha - 1/2}(x + y)dy,\qquad x\neq 0,\\
 V_\alpha(f)(0)=f(0).
 \end{gather*}

\begin{remark}
 For $\alpha
>-1/2$, we have
\[E_\alpha(\lambda.)=V_\alpha(e^{\lambda.}),\qquad \lambda\in \mathbb{C}.\]
\end{remark}

 \begin{proposition}[see \cite{Trimeche2}, Theorem 6.3]\label{prop1}
The operator $V_\alpha$ is a topological automorphism of~${\cal
E}(\mathbb{R})$, and satisfies the transmutation relation:
\[\Lambda_\alpha(V_\alpha(f)) = V_\alpha\left(\frac{d}{dx} f\right),\qquad f \in
{\cal E}(\mathbb{R}).\]
\end{proposition}

Let $\alpha >-1/2$ and we
def\/ine the dual Dunkl intertwining operator $^tV_\alpha$ on ${\cal
S}(\mathbb{R})$ (the Schwartz space on $\mathbb{R}$), by
\[ ^tV_\alpha(f)(x):= a_\alpha\int_{|y|\geq |x| }\mbox{sgn}
(y)\big(y^2-x^2\big)^{\alpha- 1/2}(x+y)f(y)dy,\]
 which can be written as:
\[ ^tV_\alpha(f)(x)= a_\alpha\,\mbox{sgn}
(x)|x|^{2\alpha+1}\int_{|t|\geq 1}\mbox{sgn} (t)\big(t^2-1\big)^{\alpha-
1/2}(1+t)f(xt)dt.\]

 \begin{proposition}[see \cite{Trimeche3}, Theorems 3.2, 3.3]\label{prop2} \qquad {}

$(i)$ The operator $^tV_\alpha$ is a topological automorphism of
${\cal S}(\mathbb{R})$, and satisfies the transmutation
relation:\[{}^tV_\alpha(\Lambda_\alpha f) =
\frac{d}{dx}({}^tV_\alpha (f)),\qquad f\in {\cal S}(\mathbb{R}).\]

$(ii)$ For all $f \in {\cal E}(\mathbb{R})$ and $g \in {\cal
S}(\mathbb{R})$, we have
\[\int_{\mathbb{R}}V_\alpha(f)(x)g(x)|x|^{2\alpha+1}dx=\int_{\mathbb{R}}f(x)
\,{}^tV_\alpha(g)(x)dx.
\]
\end{proposition}

\begin{remark}[see \cite{Soltani}] \qquad {}

 $(i)$ For $\alpha
>-1/2$ and $f\in {\cal E}(\mathbb{R})$, we can write
\[V_\alpha(f)(x)=\Re_\alpha(f_e)(|x|)+\frac{1}{x}\Re_\alpha(Mf_o)(|x|),
\]
where
 \[f_e(x)=\frac{1}{2}(f(x)+f(-x)),\qquad
f_o(x)=\frac{1}{2}(f(x)-f(-x)),\qquad Mf_o(x)=xf_o(x),
\] and $\Re_\alpha$ is the
Riemann--Liouville transform (see \cite{Trimeche1}, page~75) given by
\[
\Re_\alpha(f_e)(x):=2a_\alpha\int^1_0f_e(xt)\big(1-t^2\big)^{\alpha-1/2}dt,\qquad x\geq 0.\]

Thus, we obtain
\[
V^{-1}_\alpha(f)(x)=\Re^{-1}_\alpha(f_e)(|x|)+\frac{1}{x}\Re^{-1}_\alpha(Mf_o)(|x|).
\]
Therefore (see also \cite{Xu}, Proposition~2.2), we get
\begin{gather*}
V^{-1}_\alpha(f_e)(x)=d_\alpha\frac{d}{dx}\left(\frac{d}{xdx}\right)^{r}\left\{x^{2r+1}
\int^1_0f_e(xt)\big(1-t^2\big)^{r-\alpha-1/2}t^{2\alpha+1}dt\right\},\\
V^{-1}_\alpha(f_o)(x)=d_\alpha\left(\frac{d}{xdx}\right)^{r+1}\left\{x^{2r+2}
\int^1_0f_o(xt)\big(1-t^2\big)^{r-\alpha-1/2}t^{2\alpha+2}dt\right\},
\end{gather*}
where $r=[\alpha+1/2]$ denote the integer part of $\alpha+1/2$, and
$d_\alpha=\frac{2^{-r}\pi}{\Gamma(\alpha+1)\Gamma(r-\alpha+1/2)}$.

$(ii)$ For $\alpha
>-1/2$ and $f\in {\cal S}(\mathbb{R})$, we can write \[{}^tV_\alpha(f)(x)=W_\alpha(f_e)(|x|)+xW_\alpha(M^{-1}f_o)(|x|),
\]
where
\[
M^{-1}f_o(x)=\frac{1}{2x}(f(x)-f(-x)),\] and $W_\alpha$ is
the Weyl integral transform (see \cite[page~85]{Trimeche1}) given by
\[W_\alpha(f_e)(x):=2a_\alpha x^{2\alpha+1}\int^\infty_1f_e(xt)\big(t^2-1\big)^{\alpha-1/2}tdt,\qquad x\geq 0.\]
Thus, we obtain
\[({}^tV_\alpha)^{-1}f(x)=W^{-1}_\alpha(f_e)(|x|)+xW^{-1}_\alpha(M^{-1}f_o)(|x|).\]
\end{remark}

The Dunkl kernel gives rise to an integral transform, called Dunkl
transform on $\mathbb{R}$, which was introduced by Dunkl in~\cite{Dunkl3},
where already many basic properties were established. Dunkl's
results were completed and extended later on by de Jeu in~\cite{deJeu}.

 The
Dunkl transform of a function $f \in {\cal S}(\mathbb{R})$, is given
by
\[{\cal F}_\alpha(f)(\lambda) := \int_{\mathbb{R}} E_\alpha (- i\lambda x) f(x)
|x|^{2\alpha+1}dx, \qquad \lambda\in \mathbb{R}.\]
 We notice that ${\cal
F}_{-1/2}\,$ agrees with the Fourier transform ${\cal F}$ that is given by:
\[{\cal F}(f)(\lambda) : = \int_{\mathbb{R}} e^{- i\lambda x} f(x)\,
dx, \qquad \lambda \in \mathbb{R}.\]

 \begin{proposition}[see \cite{deJeu}]\label{prop3} \qquad {}

$(i)$ For all $f \in {\cal S}(\mathbb{R})$, we have \[{\cal
F}_\alpha(\Lambda_\alpha f)(\lambda) = i\lambda\,{\cal F}_\alpha
(f)(\lambda),\qquad \lambda \in \mathbb{R},\] where $\Lambda_\alpha$
is the Dunkl operator given by \eqref{eq1}.

$(ii)$ ${\cal F}_\alpha$ possesses on ${\cal S}(\mathbb{R})$ the
following decomposition: \[{\cal F}_\alpha (f) = {\cal
F}\,\circ\,^tV_\alpha(f),\qquad f \in {\cal S}(\mathbb{R}).\]

$(iii)$ ${\cal F}_\alpha$ is a topological automorphism of ${\cal
S}(\mathbb{R})$, and for $f\in {\cal S}(\mathbb{R})$ we have
\[f(x)=c_\alpha\int_{\mathbb{R}} E_\alpha (i\lambda
x) {\cal F}_\alpha(f)(\lambda) |\lambda|^{2\alpha+1}d\lambda,
\]
where
\[
c_\alpha=\frac{1}{[2^{\alpha+1}\Gamma(\alpha+1)]^2}.\]

$(iv)$ The normalized Dunkl transform $\sqrt{c_\alpha}\,{\cal
F}_\alpha$ extends uniquely to an isometric isomorphism of
$L^2(\mathbb{R},|x|^{2\alpha+1}dx)$ onto itself. In particular,
\[\int_{\mathbb{R}}
|f(x)|^2 |x|^{2\alpha+1}dx=c_\alpha\int_{\mathbb{R}} |{\cal
F}_\alpha(f)(\lambda)|^2 |\lambda|^{2\alpha+1}d\lambda.\]
\end{proposition}

For $T\in {\cal S}'(\mathbb{R})$, we def\/ine the Dunkl transform
${\cal F}_\alpha(T)$ of $T$, by
\begin{equation}
\langle{\cal F}_\alpha(T),
\varphi\rangle : =\langle T, {\cal F}_\alpha(\varphi)\rangle,\qquad
\varphi \in {\cal S}(\mathbb{R}).\label{eq3}
\end{equation} Thus the transform
${\cal F}_\alpha$
 extends to a topological automorphism on ${\cal S}'(\mathbb{R})$.

 \medskip

 In \cite{Trimeche3}, the author def\/ines:

 \medskip

 $\bullet\;\,$ The Dunkl translation operators $\tau_x$, $x\in \mathbb{R}$, on
${\cal E}(\mathbb{R})$, by
\[
\tau_x f(y) : =
(V_\alpha)_x\otimes(V_\alpha)_y\big[(V_\alpha)^{-1}(f)(x+y)\big],\qquad y\in
\mathbb{R}.
\] These operators satisfy for $x, y \in \mathbb{R}$ and
$\lambda \in \mathbb{C}$ the following properties:
\begin{gather*}
E_\alpha(\lambda x) E_\alpha(\lambda y) = \tau_x(E_\alpha(\lambda .))(y), \qquad \mbox{and}\\
 {\cal F}_\alpha(\tau_x f)(\lambda) = E_k(i\lambda x) {\cal F}_\alpha (f)(\lambda),\qquad f
\in {\cal S}(\mathbb{R}).
\end{gather*}

\begin{proposition}[see~\cite{Rosler1}] \label{prop4} If
$f\in {\cal C}(\mathbb{R})$ (the space of continuous functions on
$\mathbb{R}$) and $x,y \in \mathbb{R}$ such that $(x,y)\neq (0,0)$,
then
\begin{gather*}
\tau_x f(y)= a_\alpha\int^\pi_{0}
\left[f_e((x,y)_\theta)+f_o((x,y)_\theta)\frac{x+y}{(x,y)_\theta}\right][1-\mbox{\rm sgn}(xy)\cos
\theta] \sin^{2\alpha}\theta d\theta,\\
f_e(z)=\tfrac{1}{2}(f(z)+f(-z)),\qquad f_o(z)=\tfrac{1}{2}(f(z)-f(-z)),\\
(x,y)_\theta=\sqrt{x^2+y^2-2|xy|\cos \theta}.
\end{gather*}
\end{proposition}

$\bullet\;\,$ The Dunkl convolution product $\ast_\alpha$ of two
functions $f$ and $g$ in ${\cal S}(\mathbb{R})$, by
 \[f \ast_\alpha
g(x) : = \int_{\mathbb{R}} \tau_x f(-y)
g(y)|y|^{2\alpha+1}dy ,\qquad x \in \mathbb{R}.\]
 This convolution
is associative, commutative in ${\cal S}(\mathbb{R})$ and satisf\/ies (see \cite[Theorem~7.2]{Trimeche3}):
\[{\cal F}_\alpha(f \ast_\alpha g) = {\cal F}_\alpha (f) {\cal
F}_\alpha(g).\]

For $T \in{\cal S}'(\mathbb{R})$ and $f \in {\cal S}(\mathbb{R})$,
we def\/ine the Dunkl convolution product $T \ast_\alpha f$, by
\begin{equation} T
\ast_\alpha f (x) : = \langle T (y), \tau_x f (-y)\rangle,\qquad x
\in \mathbb{R}.\label{eq4}
\end{equation}

Note that $\ast_{-1/2}$ agrees with the standard convolution
$\ast$:
\[T \ast f (x): = \langle T(y), f(x-y)\rangle.\]

\section{The Dunkl Sonine transform}\label{section3}

In this section we study the Dunkl Sonine transform, which also studied by Y.~Xu
on polynomials in~\cite{Xu}. For thus we consider the following identity, which
is a consequence of Xu's result when we extend the result of Lemma~2.1 on ${\cal E}(\mathbb{R})$.
\begin{proposition}\label{prop5} Let $\alpha,  \beta \in\;  ]{-}1/2, \infty[$,
 such that $\beta>\alpha$. Then
 \begin{equation} E_\beta(\lambda x)
 =a_{\alpha,\beta}
\int^{1}_{-1}E_\alpha(\lambda xt)\big(1 - t^2\big)^{\beta-\alpha - 1}(1 +
t)|t|^{2\alpha+1}dt,\label{eq5}
\end{equation}
where
\[a_{\alpha,\beta}=\frac{\Gamma(\beta+1)}{\Gamma(\beta-\alpha)\Gamma(\alpha+1)}.\]
\end{proposition}

\begin{proof} From \eqref{eq2}, we have
\[
\int^{1}_{-1}E_\alpha(\lambda xt)\big(1 - t^2\big)^{\beta-\alpha - 1}(1 +
t)|t|^{2\alpha+1}dt=\sum^{\infty}_{n=0} \displaystyle{\frac{(\lambda
x)^n}{b_n(\alpha)}}\;I_n(\alpha,\beta) ,
\]
where \[I_n(\alpha,\beta)=\int^{1}_{-1}t^n\big(1
- t^2\big)^{\beta-\alpha - 1}(1 + t)|t|^{2\alpha+1}dt,\]
or
\begin{gather*}I_{2n}(\alpha,\beta)=2\int^{1}_{0}\big(1
- t^2\big)^{\beta-\alpha - 1}t^{2n+2\alpha+1}dt =\int^{1}_{0}(1
- y)^{\beta-\alpha - 1}y^{n+\alpha}dy \\
\phantom{I_{2n}(\alpha,\beta)}{} =\frac{\Gamma(\beta-\alpha)\Gamma(n+\alpha+1)}{\Gamma(n+\beta+1)},
\end{gather*}and
\[
I_{2n+1}(\alpha,\beta)=2\int^{1}_{0}\big(1
- t^2\big)^{\beta-\alpha - 1}t^{2n+2\alpha+3}dt=I_{2n}(\alpha+1,\beta+1) .
\]
Thus
\[\int^{1}_{-1}E_\alpha(\lambda xt)\big(1 - t^2\big)^{\beta-\alpha - 1}(1 +
t)|t|^{2\alpha+1}dt=\frac{\Gamma(\beta-\alpha)\Gamma(\alpha+1)}{\Gamma(\beta+1)}\,E_\beta(\lambda
x),\] which gives the desired result.
\end{proof}

\begin{remark} We can write the formula \eqref{eq5} by the following
\[
E_\beta(\lambda x)
 =a_{\alpha,\beta}\;
\mbox{sgn}(x)\,|x|^{-(2\beta+1)}\int^{|x|}_{-|x|}E_\alpha(\lambda
y)\big(x^2 - y^2\big)^{\beta-\alpha - 1}(x + y)|y|^{2\alpha+1}dy,\qquad x\neq 0.\]
\end{remark}

\begin{definition} \label{def1} Let $\alpha, \beta \in \;]{-}1/2, \infty[$,
 such that $\beta>\alpha$. We def\/ine the Dunkl Sonine trans\-form~$S_{\alpha, \beta}$ on
${\cal E}(\mathbb{R})$, by
\[S_{\alpha, \beta}(f)(x)
:=a_{\alpha,\beta} \int^{1}_{-1}f(xt)\big(1 - t^2\big)^{\beta-\alpha - 1}(1
+ t)|t|^{2\alpha+1}dt,\] which can be written as:
\begin{gather*}
S_{\alpha, \beta}(f)(x)
=a_{\alpha,\beta}\;
\mbox{sgn}(x)\,|x|^{-(2\beta+1)}\int^{|x|}_{-|x|}f(y)\big(x^2 -
y^2\big)^{\beta-\alpha - 1}(x + y)|y|^{2\alpha+1}dy,\qquad x\neq 0,\\
 S_{\alpha, \beta}(f)(0)=f(0).
 \end{gather*}
\end{definition}

\begin{remark} For $\alpha, \beta \in \;]{-}1/2, \infty[$,
 such that $\beta>\alpha$, we have
  \begin{equation}
 E_\beta(\lambda.)=S_{\alpha, \beta}(E_\alpha(\lambda.)),
 \qquad \lambda\in \mathbb{C}.\label{eq6}
 \end{equation}
 \end{remark}

 \begin{definition}\label{def2} Let $\alpha, \beta \in
\;]{-}1/2, \infty[$, such that $\beta >\alpha$. We def\/ine the dual
Dunkl Sonine transform $^tS_{\alpha, \beta}$ on ${\cal
S}(\mathbb{R})$, by
\[ ^tS_{\alpha, \beta}(f)(x):= a_{\alpha,\beta}\int_{|y|\geq |x|}\mbox{sgn}
(y)\big(y^2-x^2\big)^{\beta-\alpha- 1}(x+y)f(y)dy,\] which can be written
as:
\[ ^tS_{\alpha, \beta}(f)(x)= a_{\alpha,\beta}\,\mbox{sgn}
(x)|x|^{2(\beta-\alpha)}\int_{|t|\geq 1}\mbox{sgn}
(t)\big(t^2-1\big)^{\beta-\alpha- 1}(t+1)f(xt)dt.\]
\end{definition}

\begin{proposition}\label{prop6} \qquad {}

$(i)$ For all $f \in {\cal E}(\mathbb{R})$ and $g \in {\cal
S}(\mathbb{R})$, we have
\[\int_{\mathbb{R}}S_{\alpha,\beta}(f)(x)g(x)|x|^{2\beta+1}dx=\int_{\mathbb{R}}f(x)
\,{}^tS_{\alpha,\beta}(g)(x)|x|^{2\alpha+1}dx.\]

$(ii)$ ${\cal F}_\beta$ possesses on ${\cal S}(\mathbb{R})$ the
following decomposition: \[{\cal F}_\beta (f) = {\cal
F}_\alpha\,\circ\,^tS_{\alpha,\beta}(f),\qquad f \in {\cal
S}(\mathbb{R}).\]
\end{proposition}

\begin{proof}
 Part $(i)$ follows from Def\/inition~\ref{def1}
by Fubini's theorem. Then part $(ii)$ follows from $(i)$ and \eqref{eq6} by
taking $f=E_\alpha(-i\lambda.)$.\end{proof}

In  \cite[Lemma 2.1]{Xu} Y.~Xu proves the identity $S_{\alpha, \beta}= V_\beta\,\circ\,V^{-1}_\alpha$
 on polynomials. As the intertwiner is a homeomorphism on ${\cal E}(\mathbb{R})$
  and polynomials are dense in ${\cal E}(\mathbb{R})$, this gives the identity also
on ${\cal E}(\mathbb{R})$. In the following we give a second method to prove this identity.

\begin{theorem}\label{th1}\qquad{}

$(i)$ The operator $^tS_{\alpha, \beta}$ is a topological
automorphism of ${\cal S}(\mathbb{R})$, and satisfies the following
relations:
\begin{gather*}
{}^tS_{\alpha, \beta}(f) = (^tV_\alpha)^{-1}\,\circ\,^tV_\beta(f),\qquad f \in
{\cal S}(\mathbb{R}),\\
  ^tS_{\alpha, \beta}(\Lambda_\beta f)=\Lambda_\alpha(^tS_{\alpha, \beta}(f)),\qquad f \in
{\cal S}(\mathbb{R}).
\end{gather*}

 $(ii)$ The operator
$S_{\alpha, \beta}$ is a topological automorphism of ${\cal
E}(\mathbb{R})$, and satisfies the following relations:
\begin{gather*}
S_{\alpha, \beta}(f) = V_\beta\,\circ\,V^{-1}_\alpha(f),\qquad f \in
{\cal E}(\mathbb{R}),\\
 \Lambda_\beta(S_{\alpha, \beta}(f)) = S_{\alpha, \beta}(\Lambda_\alpha f),\qquad f \in
{\cal E}(\mathbb{R}).
\end{gather*}
\end{theorem}

\begin{proof}
$(i)$ From Proposition~\ref{prop6} $(ii)$, we have
\begin{equation} ^tS_{\alpha,\beta}(f)=({\cal
F}_\alpha)^{-1}\,\circ\,{\cal F}_\beta (f).\label{eq7}
\end{equation}
 Using
Proposition~\ref{prop3} $(ii)$, we obtain
\begin{equation}
^tS_{\alpha, \beta}(f) =
(^tV_\alpha)^{-1}\,\circ\,^tV_\beta(f),\qquad f \in {\cal
S}(\mathbb{R}).\label{eq8}
\end{equation}
 Thus from Proposition~\ref{prop2} $(i)$,
\[ ^tS_{\alpha, \beta}(\Lambda_\beta f) =
(^tV_\alpha)^{-1}\,\circ\,^tV_\beta(\Lambda_\beta f)=
(^tV_\alpha)^{-1}\left(\frac{d}{dx}\,^tV_\beta(f)\right).\] Using the fact that
\[{}^tV_\alpha(\Lambda_\alpha f) =
\frac{d}{dx}({}^tV_\alpha (f))\ \Longleftrightarrow \ \Lambda_\alpha
({}^tV_\alpha)^{-1}(f) = ({}^tV_\alpha)^{-1}\left(\frac{d}{dx}f\right),\] we
obtain
\[ ^tS_{\alpha, \beta}(\Lambda_\beta f) =
\Lambda_\alpha(^tV_\alpha)^{-1}(\,^tV_\beta(f))=\Lambda_\alpha(^tS_{\alpha,
\beta}(f)).\]

$(ii)$ From Proposition~\ref{prop2} $(ii)$, we have
\[\int_{\mathbb{R}}f(x)
\,{}^tV_\beta(g)(x)dx=\int_{\mathbb{R}}V_\beta(f)(x)g(x)|x|^{2\beta+1}dx.\]
On other hand, from \eqref{eq8}, Proposition~\ref{prop2} $(ii)$ and Proposition~\ref{prop6} $(i)$
we have
\begin{gather*}
\int_{\mathbb{R}}f(x)
\,{}^tV_\beta(g)(x)dx=\int_{\mathbb{R}}f(x)\,{}^tV_\alpha\,\circ\,^tS_{\alpha,\beta}(g)(x)dx
=\int_{\mathbb{R}}V_\alpha(f)(x)\,^tS_{\alpha,\beta}(g)(x)|x|^{2\alpha+1}dx
\\
\phantom{\int_{\mathbb{R}}f(x)
\,{}^tV_\beta(g)(x)dx}{} =\int_{\mathbb{R}}S_{\alpha,\beta}\,\circ\,V_\alpha(f)(x)g(x)
|x|^{2\beta+1}dx.
\end{gather*}
Then \[S_{\alpha,\beta}\,\circ\,V_\alpha(f)=V_\beta(f).\] Hence from
Proposition~\ref{prop1}, \[\Lambda_\beta(S_{\alpha, \beta}(f)) = \Lambda_\beta
V_\beta(V^{-1}_\alpha(f)) = V_\beta\left(\frac{d}{dx}V^{-1}_\alpha(f)\right).\]
Using the fact that
\[\Lambda_\alpha(V_\alpha(f)) =
V_\alpha\left(\frac{d}{dx} f\right) \ \Longleftrightarrow \ V^{-1}_\alpha(\Lambda_\alpha
f)= \frac{d}{dx}V^{-1}_\alpha(f),\]
we obtain
\[\Lambda_\beta(S_{\alpha, \beta}(f)) = V_\beta\,\circ\,V^{-1}_\alpha(\Lambda_\alpha
f)=S_{\alpha, \beta}(\Lambda_\alpha f),\] which completes the proof
of the theorem.\end{proof}

 \section[Complex powers of $\Delta_\alpha$]{Complex powers of $\boldsymbol{\Delta_\alpha}$}\label{section4}

For $\lambda\in \mathbb{C}$, $\mbox{Re}(\lambda)>-1$, we denote by
$|x|^{\lambda}$ the tempered distribution def\/ined by
\begin{equation}
\langle |x|^{\lambda}, \varphi \rangle : =
\int_{\mathbb{R}}|x|^{\lambda}\varphi(x)dx,\qquad \varphi \in {\cal
S}(\mathbb{R}).\label{eq9}
\end{equation}
We write \[\langle |x|^{\lambda}, \varphi \rangle  =
\int^\infty_{0}x^{\lambda}[\varphi(x)+\varphi(-x)]dx,\qquad \varphi \in {\cal
S}(\mathbb{R}),\] then from \cite{Rachdi}, we obtain the following result.

\begin{lemma}\label{lem1} Let $\varphi \in {\cal S}(\mathbb{R})$. The mapping
$g:\lambda \rightarrow \langle|x|^{\lambda},\varphi\rangle$ is
complex-valued function and has an analytic extension to
$\mathbb{C}\backslash \{-(1+2\ell),\;\ell \in
\mathbb{N}\}$, with simple poles $-(2\ell+1)$, $\ell \in
\mathbb{N}$ and \[\mbox{\rm Res}(g,-1-2\ell)=2\frac{\varphi^{(2\ell)}(0)}{(2\ell)!}.\]
\end{lemma}

\begin{proposition}\label{prop7} Let $\varphi \in {\cal S}(\mathbb{R})$.

$(i)$ The function $\lambda \rightarrow \langle
|x|^{\lambda+2\alpha+1},\varphi \rangle$ is analytic on
$\mathbb{C}\backslash \{-(2\alpha+2\ell+2),\;\ell \in
\mathbb{N}\}$, with simple poles $-(2\alpha+2\ell+2)$, $\ell \in
\mathbb{N}$.

$(ii)$ The function $\lambda \rightarrow \frac{ 2^{2\alpha +
\lambda+2}\Gamma(\alpha+1)\Gamma(\frac{2\alpha+\lambda+2}{2})}{\Gamma
(-\lambda/2)} \langle |x|^{-(\lambda+1)},\varphi\rangle$ is analytic
on $\mathbb{C}\backslash \{-(2\alpha+2\ell+2),$ $\ell \in
\mathbb{N}\}$, with simple poles $-(2\alpha+2\ell+2)$, $\ell \in
\mathbb{N}$.

 $(iii)$ For $\lambda \in \mathbb{C}\backslash \{-(2\alpha+2\ell+2),\;\ell \in
\mathbb{N}\}$ we have \[{\cal
F}_\alpha\big(|x|^{\lambda+2\alpha+1}\big)=\frac{ 2^{2\alpha +
\lambda+2}\Gamma(\alpha+1)\Gamma(\frac{2\alpha+\lambda+2}{2})}{\Gamma
(-\lambda/2)}|x|^{-(\lambda+1)},\qquad \mbox{in ${\cal S}'$-sense}.\]

$(iv)$ For $\lambda \in \mathbb{C}\backslash
\{-(2\alpha+2\ell+2),\;\ell \in \mathbb{N}\}$ we have
\[|x|^{\lambda+2\alpha+1}=\frac{ 2^{
\lambda}\,\Gamma(\frac{2\alpha +
\lambda+2}{2})}{\Gamma(\alpha+1)\Gamma (-\lambda/2)}{\cal
F}_\alpha(|x|^{-(\lambda+1)}),\qquad \mbox{in ${\cal S}'$-sense}.\]\end{proposition}

\begin{proof}

$(i)$ Follows directly from Lemma~\ref{lem1}.

$(ii)$ From  \cite[pages~2 and 8]{Lebedev}  the function $\lambda \rightarrow \Gamma
(\frac{2\alpha+\lambda + 2}{2})$ has an analytic extension to
$\mathbb{C}\backslash \{-(2\alpha+2\ell+2),\;\ell \in \mathbb{N}\}$, with simple poles
$-(2\alpha+2\ell+2)$, $\ell \in \mathbb{N}$,
and the function
$\lambda \rightarrow \frac{1}{\Gamma (-\lambda/2)}$ has zeros
$2\ell$, $\ell \in \mathbb{N}$. Thus from Lemma~\ref{lem1} we see that
\[\lambda \rightarrow \frac{ 2^{2\alpha +
\lambda+2}\Gamma(\alpha+1)\Gamma(\frac{2\alpha+\lambda+2}{2})}{\Gamma
(-\lambda/2)} \langle |x|^{-(\lambda+1)},\varphi\rangle\] is
analytic on $\mathbb{C}\backslash
\{-(2\alpha+2\ell+2),\;\ell \in \mathbb{N}\}$, with simple poles
$-(2\alpha+2\ell+2)$, $\ell \in \mathbb{N}$.\vspace{1mm}

 $(iii)$ Let determine the value of ${\cal F}_\alpha(
|x|^{\lambda+2\alpha+1})$ in the ${\cal S}'$-sense. We put
$\psi_t(x):=e^{-tx^2}$, $t>0$. Then $\psi_t \in {\cal
S}(\mathbb{R})$, and from~\cite{Rosler2}:
\[{\cal
F}_\alpha(\psi_t)(x)=\Gamma(\alpha+1)t^{-(\alpha+1)}e^{-x^2/4t},\qquad
x \in \mathbb{R}.\]

Furthermore, for $\varphi \in{\cal S}(\mathbb{R})$ we have
\[\int_{\mathbb{R}}{\cal F}_\alpha(\varphi)(x)\psi_t(x)|x|^{2\alpha+1}dx =
\Gamma(\alpha+1) \int_{\mathbb{R}}\varphi(x)
t^{-(\alpha+1)}e^{-x^2/4t}|x|^{2\alpha+1}dx.\]
Multiplying both sides by $t^{-\lambda/2-1}$ and integrating over
$(0, \infty)$, we obtain for $\mbox{Re}(\lambda)\in \,]{-}(2\alpha + 2),
0[$:
\[\int_{\mathbb{R}}{\cal
F}_\alpha(\varphi)(x)|x|^{\lambda+2\alpha+1}dx = \frac{ 2^{2\alpha +
\lambda+2}\Gamma(\alpha+1)\Gamma(\frac{2\alpha+\lambda+2}{2})}{\Gamma
(-\lambda/2)}\int_{\mathbb{R}}\varphi (x)|x|^{-(\lambda+1)}dx.\]
This and from \eqref{eq3} we get for $\mbox{Re}(\lambda)\in \,]{-}(2\alpha + 2), 0[$:
\[{\cal F}_\alpha
(|x|^{\lambda+2\alpha+1})=\frac{ 2^{2\alpha +
\lambda+2}\Gamma(\alpha+1)\Gamma(\frac{2\alpha+\lambda+2}{2})}{\Gamma
(-\lambda/2)}|x|^{-(\lambda+1)}.\] The result follows by analytic
continuation.

$(iv)$ From $(iii)$ we have
\[|x|^{\lambda+2\alpha+1}=\frac{ 2^{2\alpha +
\lambda+2}\Gamma(\alpha+1)\Gamma(\frac{2\alpha+\lambda+2}{2})}{\Gamma
(-\lambda/2)}{\cal F}^{-1}_\alpha\big(|x|^{-(\lambda+1)}\big).\] Using the
fact that
\[\langle{\cal F}^{-1}_\alpha\big(|x|^{-(\lambda+1)}\big),
\varphi\rangle =\langle |x|^{-(\lambda+1)}, {\cal
F}^{-1}_\alpha(\varphi)\rangle,\qquad \varphi \in {\cal
S}(\mathbb{R}).\]
 By applying \eqref{eq9} and Proposition~\ref{prop3} $(iii)$, we obtain
\[\langle{\cal F}^{-1}_\alpha(|x|^{-(\lambda+1)}),
\varphi\rangle =c_\alpha\int_{\mathbb{R}}|x|^{-(\lambda+1)}{\cal
F}_\alpha(\varphi)(-x)dx,\qquad \varphi \in {\cal S}(\mathbb{R}).\]
Then
\[{\cal F}^{-1}_\alpha\big(|x|^{-(\lambda+1)}\big)=c_\alpha\,{\cal
F}_\alpha\big(|x|^{-(\lambda+1)}\big),\] which gives the
result.\end{proof}

\begin{definition}\label{def3} For $\lambda\in \mathbb{C}\backslash
\{-(\alpha+\ell+1),\;\ell \in \mathbb{N}\}$, the
complex powers of the Dunkl Laplacian~$\Delta_\alpha$ are def\/ined
for $f \in {\cal S}(\mathbb{R})$ by \[(-\Delta_\alpha)^{\lambda}f(x)
: =
\frac{2^{2\lambda}\Gamma(\alpha+\lambda+1)}{\Gamma(\alpha+1)\Gamma
(-\lambda)}|x|^{-(2\lambda+1)}\ast_\alpha f(x),\]where $\ast_\alpha$
is the Dunkl convolution product given by \eqref{eq4}.\end{definition}

In the next part of this section we use Def\/inition~\ref{def3} and Proposition~\ref{prop7} $(iv)$
 to establish the following result:
 \[{\cal
F}_\alpha\big((-\Delta_\alpha)^{\lambda}f\big)(x)= |x|^{2\lambda}{\cal
F}_\alpha(f)(x).
\]
\begin{proposition}\label{prop8}
 For $\lambda\in \mathbb{C}\backslash
\{-(\alpha+\ell+1),\;\ell \in \mathbb{N}\}$ and $f \in
{\cal S}(\mathbb{R})$, \[(-\Delta_\alpha)^{\lambda}f(x) =
b_\alpha(\lambda)\int_{\mathbb{R}}\left[\int^\pi_{0}
\frac{(1+\mbox{\rm sgn}(xy)\cos
\theta)}{(x,y)^{2(\lambda+\alpha+1)}_\theta} \sin^{2\alpha}\theta
d\theta\right]f(y)|y|^{2\alpha+1}dy,\]
where
\[b_\alpha(\lambda)=\frac{2^{2\lambda}\Gamma(\alpha+\lambda+1)}{\sqrt{\pi}\,\Gamma(\alpha+1/2)\Gamma
(-\lambda)} ,\qquad (x,y)_\theta=\sqrt{x^2+y^2-2|xy|\cos \theta}.\]
\end{proposition}

\begin{proof}
 From Def\/inition~\ref{def3}, \eqref{eq4} and \eqref{eq9}, we have
\begin{gather*}
(-\Delta_\alpha)^{\lambda}f(x)
=\frac{2^{2\lambda}\Gamma(\alpha+\lambda+1)}{\Gamma(\alpha+1)\Gamma
(-\lambda)}\langle|y|^{-(2\lambda+1)},\tau_xf(-y) \rangle \\
\phantom{(-\Delta_\alpha)^{\lambda}f(x)}{} =\frac{2^{2\lambda}\Gamma(\alpha+\lambda+1)}{\Gamma(\alpha+1)\Gamma
(-\lambda)}\int_{\mathbb{R}}|y|^{-2(\lambda+\alpha+1)}\tau_xf(-y)|y|^{2\alpha+1}dy.
\end{gather*} So
\[(-\Delta_\alpha)^{\lambda}f(x) =\int_{\mathbb{R}}\tau_x(|y|^{-2(\lambda+\alpha+1)})(-y)f(y)|y|^{2\alpha+1}dy.\]
Then the result follows from Proposition~\ref{prop4}.\end{proof}

\begin{note} We denote by

$\bullet\;\,\Psi$ the subspace of ${\cal S}(\mathbb{R})$ consisting
of functions $f$, such that
\[f^{(k)}(0)=0,\qquad \forall\;k \in \mathbb{N}.\]

$\bullet\;\, \Phi_\alpha$ the subspace of ${\cal S}(\mathbb{R})$
consisting of functions $f$, such that
\[\int_{\mathbb{R}}f(y)\,y^k|y|^{2\alpha+1}dy=0,\qquad
\forall\;k\in \mathbb{N}.\] The spaces $\Psi$ and $\Phi_{-1/2}$ are
well-known in the literature as Lizorkin spaces (see \cite{Rachdi,Nessibi,Samko}).
\end{note}

\begin{lemma}[see \cite{Rachdi}] \label{lem2} The multiplication operator $M_\lambda: f \rightarrow
|x|^\lambda f$, $\lambda \in \mathbb{C}$, is a topological
automorphism of $\Psi$. Its inverse operator is $(M_\lambda)^{-1} =
M_{-\lambda}$.
\end{lemma}

\begin{theorem}\label{th2}\qquad{}

$(i)$ The Dunkl transform ${\cal F}_\alpha$ is a topological
isomorphism from $\Phi_\alpha$ onto $\Psi$.

$(ii)$ The operator ${}^tS_{\alpha,\beta}$ is a topological
isomorphism from $\Phi_\beta$ onto $\Phi_\alpha$.

$(iii)$ For $\lambda\in \mathbb{C}\backslash
\{-(\alpha+\ell+1), \;\ell \in \mathbb{N}\}$ and $f\in
\Phi_\alpha\,$, the function $(-\Delta_\alpha)^{\lambda}f$ belongs
to $\in \Phi_\alpha$, and
\begin{equation} {\cal
F}_\alpha((-\Delta_\alpha)^{\lambda}f)(x)= |x|^{2\lambda}{\cal
F}_\alpha(f)(x).\label{eq10}
\end{equation}
\end{theorem}

\begin{proof}

$(i)$ Let $f \in \Phi_\alpha$, then \[({\cal F}_\alpha(f))^{(k)}(0) =
(-i)^{k}\frac{k!}{b_k(\alpha)}\int_{\mathbb{R}}f(x)\,x^k|x|^{2\alpha+1}dy=0,\qquad
\forall\; k \in \mathbb{N}.\] Hence ${\cal F}_\alpha(f)\in \Psi$.

 Conversely, let $g \in \Psi$. Since ${\cal F}_\alpha$ is a topological automorphism
of ${\cal S}(\mathbb{R})$. There exists $f \in {\cal
S}(\mathbb{R})$, such that ${\cal F}_\alpha(f)=g$. Thus
\[g^{(k)}(0) =
(-i)^{k}\frac{k!}{b_k(\alpha)}\int_{\mathbb{R}}f(x)\,x^k|x|^{2\alpha+1}dy=0,\qquad
\forall\; k \in \mathbb{N}.\] So $f \in \Phi_\alpha$ and ${\cal
F}_\alpha(f)=g$.

$(ii)$ follows directly from $(i)$ and \eqref{eq7}.

$(iii)$ Similarly to the standard convolution if $f\in {\cal
S}(\mathbb{R})$ and $S\in {\cal S}'(\mathbb{R})$, then $S
\ast_\alpha f \in {\cal E}(\mathbb{R})$ and
$T_{|x|^{2\alpha+1}\,S\ast_\alpha f} \in {\cal S}'(\mathbb{R})$.
Moreover \[{\cal F}_\alpha(T_{|x|^{2\alpha+1}\,S\ast_\alpha f}) =
{\cal F}_\alpha(f){\cal F}_\alpha(S).\]

Let $f\in \Phi_\alpha$ and $\lambda\in \mathbb{C}\backslash
\{-(\alpha+\ell+1),\;\ell \in \mathbb{N}\}$.
Consequently, from Def\/inition~\ref{def3}, Proposition~\ref{prop7}~$(iv)$ and \eqref{eq9} we have
\begin{equation}
{\cal F}_\alpha(T_{|x|^{2\alpha+1}(-\Delta_\alpha)^{\lambda}f}) =
|x|^{2\lambda+2\alpha+1}{\cal F}_\alpha(f) =
T_{|x|^{2\lambda+2\alpha+1}{\cal F}_\alpha(f)}.\label{eq11}
\end{equation}
 On the other hand from \eqref{eq3},
\begin{equation}
{\cal F}_\alpha(T_{|x|^{2\alpha+1}(-\Delta_\alpha)^{\lambda}f}) =
T_{|x|^{2\alpha+1}{\cal F
}_\alpha((-\Delta_\alpha)^{\lambda}f)}.\label{eq12}
\end{equation}
 From \eqref{eq11} and
\eqref{eq12}, we obtain
\[{\cal F}_\alpha((-\Delta_\alpha)^{\lambda}f)= |x|^{2\lambda}{\cal
F}_\alpha(f).\]
 Then by Lemma~\ref{lem2} and $(i)$ we deduce that
$(-\Delta_\alpha)^{\lambda}f \in \Phi_\alpha$. \end{proof}

\section[Inversion formulas for $ S_{\alpha,\beta}$ and
${}^tS_{\alpha,\beta}$]{Inversion formulas for $\boldsymbol{S_{\alpha,\beta}}$ and
$\boldsymbol{{}^tS_{\alpha,\beta}}$}\label{section5}

In this section, we establish inversion formulas for the Dunkl Sonine
transform and its dual.
\begin{definition}\label{def4} We def\/ine the operators $K_1$, $K_2$ and $K_3$, by
\begin{gather*}
K_1(f): =  \frac{c_\beta}{c_\alpha}\,{\cal F}^{-1}_\alpha\big(|\lambda|^{2(\beta-\alpha)}\,{\cal
F}_\alpha(f)\big)=\frac{c_\beta}{c_\alpha}\,(-\Delta_\alpha)^{\beta-\alpha}f,\qquad
f \in \Phi_\alpha,\\
K_2(f): = \frac{c_\beta}{c_\alpha}\,{\cal
F}^{-1}_\beta\big(|\lambda|^{2(\beta-\alpha)}\,{\cal
F}_\beta(f)\big)=\frac{c_\beta}{c_\alpha}\,(-\Delta_\beta)^{\beta-\alpha}f,\qquad
f \in \Phi_\beta,\\
K_3(f): =  \sqrt{\frac{c_\beta}{c_\alpha}}\,{\cal F}^{-1}_\alpha\big(|\lambda|^{\beta-\alpha}\,{\cal
F}_\alpha(f)\big)=\sqrt{\frac{c_\beta}{c_\alpha}}\,
(-\Delta_\alpha)^{(\beta-\alpha)/2}f,\qquad f \in \Phi_\alpha.
\end{gather*}
\end{definition}

\begin{lemma}\label{lem3} For all $g \in \Phi_\beta$, we have
\begin{equation}
K_1({}^tS_{\alpha,\beta})(g)=({}^tS_{\alpha,\beta})K_2(g).\label{eq13}$$
\end{equation}
\end{lemma}
\begin{proof}

 Let $g \in \Phi_\beta\,$. Using Proposition~\ref{prop6} $(ii)$,
\begin{equation*}
K_1({}^tS_{\alpha,\beta})(g)=\frac{c_\beta}{c_\alpha}\,{\cal
F}^{-1}_\alpha\big(|\lambda|^{2(\beta-\alpha)}\,{\cal
F}_\beta(g)\big)=({}^tS_{\alpha,\beta})K_2(g).\tag*{\qed}\end{equation*}
\renewcommand{\qed}{}
\end{proof}

\begin{theorem}\label{th3}\qquad{}

$(i)$ Inversion formulas: For all $f \in \Phi_\alpha$ and $g \in
\Phi_\beta$, we have the inversions formulas:
 \[(a)\ \  g = S_{\alpha,\beta} K_1({}^tS_{\alpha,\beta})(g),\qquad
 (b)\ \ f = ({}^tS_{\alpha,\beta})K_2S_{\alpha,\beta}(f).\]

$(ii)$ Plancherel formula: For all $f \in \Phi_\beta$ we have
\[\int_{\mathbb{R}} |f(x)|^2|x|^{2\beta+1}dx  = \int_{\mathbb{R}}
|K_3({}^tS_{\alpha,\beta}(f))(x)|^2 |x|^{2\alpha+1}dx.\]
\end{theorem}

\begin{proof}

$(i)$ Let $g \in \Phi_\beta$. From Proposition~\ref{prop3} $(iii)$, \eqref{eq6} and
Proposition~\ref{prop6} $(ii)$, we obtain
\begin{gather*} g =c_\beta\,\int_{\mathbb{R}}
S_{\alpha,\beta}(E_\alpha(i\lambda.))\,{\cal
F}_\beta(g)(\lambda)|\lambda|^{2\beta+1}d\lambda \\
\phantom{g} = c_\beta\,
S_{\alpha,\beta}\left[\int_{\mathbb{R}}E_\alpha(i\lambda.)\,{\cal
F}_\alpha\,\circ\,{}^tS_{\alpha,\beta}(g)(\lambda)|\lambda|^{2\beta+1}d\lambda\right]
\\
\phantom{g}
= \frac{ c_\beta}{c_\alpha}\,S_{\alpha,\beta}
\left[{\cal F}^{-1}_\alpha(|\lambda|^{2(\beta-\alpha)}\,{\cal F}_\alpha\,\circ\,{}^tS_{\alpha,\beta}(g))\right].
\end{gather*}
Thus \[g = S_{\alpha,\beta} K_1({}^tS_{\alpha,\beta})(g),\qquad g \in \Phi_\beta.\]
From the previous relation and \eqref{eq13}, we deduce the relation:
\[f = ({}^tS_{\alpha,\beta})K_2S_{\alpha,\beta}(f),\qquad f \in
\Phi_\alpha.\]

 $(ii)$ Let $f \in \Phi_\beta$.
From Proposition~\ref{prop3} $(iv)$ and Proposition~\ref{prop6} $(ii)$, we deduce that
\[\int_{\mathbb{R}} |f(x)|^2 |x|^{2\beta+1}dx =
c_\beta\,\int_{\mathbb{R}}\big| |\lambda|^{\beta-\alpha}{\cal
F}_\alpha({}^tS_{\alpha,\beta}(f))(\lambda)\big|^2
|\lambda|^{2\alpha+1}d\lambda.\] Thus we obtain
\[\int_{\mathbb{R}} |f(x)|^2|x|^{2\beta+1}dx =c_\alpha
\int_{\mathbb{R}} \big|{\cal
F}_\alpha\big(K_3({}^tS_{\alpha,\beta}(f))\big)(\lambda)\big|^2
|\lambda|^{2\alpha+1}d\lambda.\]
Then the result follows from this
identity by applying Proposition~\ref{prop3} $(iv)$.
\end{proof}

\begin{remark}
Let $f \in \Phi_\alpha$ and $g \in \Phi_\beta$.
 By writing $(a)$ and $(b)$ respectively
for the functions $S_{\alpha,\beta}(f)$ and
${}^tS_{\alpha,\beta}(g)$, we obtain
\[(c)\ \ f = K_1({}^tS_{\alpha,\beta})S_{\alpha,\beta}(f),\qquad (d)
\ \ g = K_2S_{\alpha,\beta}({}^tS_{\alpha,\beta})(g).\]
\end{remark}

\subsection*{Acknowledgements}

The author is very grateful
to the referees and editors for many critical comments on this paper.

\pdfbookmark[1]{References}{ref}
\LastPageEnding


\begin{thebibliography}{99}

\footnotesize\itemsep=0pt

\bibitem{Rachdi}
Baccar C., Hamadi N.B., Rachdi L.T.,
Inversion formulas for Riemann--Liouville transform and its dual associated with singular partial differential operators, {\it Int. J. Math. Math. Sci.} {\bf 2006} (2006), Art. ID 86238, 26~pages.

\bibitem{Dunkl1}
Dunkl C.F., Dif\/ferential-dif\/ference operators associated with
ref\/lections groups, {\it Trans. Amer. Math. Soc.} {\bf 311} (1989), 167--183.

\bibitem{Dunkl2}
Dunkl C.F., Integral kernels with ref\/lection group invariance, {\it Canad.
J. Math.} {\bf 43} (1991), 1213--1227.

\bibitem{Dunkl3}
Dunkl C.F., Hankel transforms associated to f\/inite ref\/lection groups,
{\it Contemp. Math.} {\bf 138} (1992), 123--138.

\bibitem{deJeu}
de Jeu M.F.E., The Dunkl transform, {\it Invent. Math.} {\bf 113} (1993), 147--162.

\bibitem{Lapointe}
Lapointe L., Vinet L., Exact operator solution of the Calogero--Sutherland model, {\it Comm. Math. Phys.} {\bf 178} (1996), 425--452, \href{http://arxiv.org/abs/q-alg/9509003}{q-alg/9509003}.

\bibitem{Lebedev}
Lebedev N.N., Special functions and their applications, Dover
Publications, Inc., New York, 1972.

\bibitem{Ludwig}
Ludwig D., The Radon transform on Euclidean space, {\it Comm. Pure. App.
Math.} {\bf 23} (1966), 49--81.

\bibitem{Nessibi}
Nessibi M.M., Rachdi L.T., Trim\`eche K., Ranges and inversion
formulas for spherical mean operator and its dual, {\it J. Math. Anal.
Appl.} {\bf 196} (1995), 861--884.

\bibitem{Rosenblum}
Rosenblum M., Generalized Hermite polynomials
and the Bose-like oscillator calculus,
in Nonselfadjoint Operators and Related Topics (Beer Sheva, 1992), {\it Oper. Theory Adv. Appl.}, Vol.~73, Birkh\"auser, Basel, 1994, 369--396, \href{http://arxiv.org/abs/math.CA/9307224}{math.CA/9307224}.

\bibitem{Rosler1}
R\"osler M., Bessel-type signed hypergroups on $\mathbb{R}$, in
Probability Measures on Groups
and Related Structures,~XI (Oberwolfach, 1994), Editors H.~Heyer and A.~Mukherjea, Oberwolfach, 1994, World Sci. Publ., River Edge, NJ, 1995, 292--304.

\bibitem{Rosler2}
R\"osler M., Generalized Hermite polynomials and the heat equation
for Dunkl operators, {\it Comm. Math. Phys.} {\bf 192} (1998), 519--542, \href{http://arxiv.org/abs/q-alg/9703006}{q-alg/9703006}.

\bibitem{Samko}
Samko S.G., Hypersingular integrals and their applications,
{\it Analytical Methods and Special Functions},
 Vol.~5, Taylor \& Francis, Ltd., London, 2002.

\bibitem{Solmon}
Solmon D.C., Asymptotic formulas for the dual Radon transform and applications, {\it Math.
Z.} {\bf 195} (1987), 321--343.

\bibitem{Soltani}
Soltani F., Trim\`eche K., The Dunkl intertwining operator and its
dual on $\mathbb{R}$ and applications, Preprint, Faculty of Sciences of Tunis,
Tunisia, 2000.

\bibitem{Stein}
Stein E.M., Singular integrals and dif\/ferentiability properties of
functions, Princeton University Press, Princeton, 1970.

\bibitem{Trimeche1}
Trim\`eche K., Transformation int\'egrale de Weyl et th\'eor\`eme de
Paley--Wiener associ\'es \`a un op\'erateur dif\/f\'erentiel singulier sur $(0,
\infty)$, {\it J. Math. Pures Appl. (9)} {\bf 60} (1981), 51--98.

\bibitem{Trimeche2}
Trim\`eche K., The Dunkl intertwining operator on spaces of functions
and distributions and integral representation of its dual, {\it Integral Transform. Spec. Funct.} {\bf 12} (2001), 349--374.

\bibitem{Trimeche3}
Trim\`eche K., Paley--Wiener theorems for the Dunkl transform and
Dunkl translation operators, {\it Integral Transform. Spec. Funct.} {\bf 13} (2002),
17--38.

\bibitem{Xu}
Xu Y., An integral formula for generalized Gegenbauer polynomials and
Jacobi polynomials, {\it Adv. in Appl. Math.} {\bf 29} (2002), 328--343.

\end{thebibliography}
\end{document}